\documentclass[11pt]{amsart}
\usepackage{amsmath, amssymb}
\usepackage{amsfonts}
\usepackage{mathrsfs}
\usepackage[arrow,matrix,curve,cmtip,ps]{xy}
 
\usepackage{amsthm}
 
\newtheorem{theorem}{Theorem}[section]
\newtheorem{lemma}[theorem]{Lemma}
\newtheorem{proposition}[theorem]{Proposition}
\newtheorem{corollary}[theorem]{Corollary}

\newtheorem*{theorem*}{Theorem}
\theoremstyle{remark}
\newtheorem{remark}[theorem]{Remark}
\newtheorem{definition}[theorem]{Definition}
\newtheorem{example}[theorem]{Example}

\numberwithin{equation}{section}
 
\newcommand{\Z}{\mathbb{Z}}
\newcommand{\N}{\mathbb{N}}
\newcommand{\R}{\mathbb{R}}
\newcommand{\C}{\mathbb{C}}

\newcommand{\F}{\mathcal{F}}
\newcommand{\G}{\mathcal{G}}

\newcommand{\gae}{\lower 2pt \hbox{$\, \buildrel {\scriptstyle >}\over {\scriptstyle
\sim}\,$}}
 
\newcommand{\lae}{\lower 2pt \hbox{$\, \buildrel {\scriptstyle <}\over {\scriptstyle
\sim}\,$}}
 
\newcommand{\MU}[1]{
\setbox0\hbox{$#1$}
\setbox1\hbox{$W$}
\ifdim\wd0>\wd1 #1^{\sim} \else \widetilde{#1} \fi
}

\begin{document}
\title[Frame theory for binary vector spaces]{Frame theory for binary vector spaces}
 
\author[B. G. Bodmann]{Bernhard G. Bodmann}
\author[M. Le]{My Le}
\author[L. Reza]{Letty Reza}
\author[M. Tobin]{Matthew Tobin}
\author[M. Tomforde]{Mark Tomforde}
 
\address{Department of Mathematics \\ University of Houston \\ Houston, TX 77204-3008 \\USA}
\email{bg{}b{}@{}math{}.uh.{}edu}
 
\address{Department of Mathematics \\ University of Houston \\ Houston, TX 77204-3008 \\USA}
\email{myle\_64@yahoo.com}
 
\address{Department of Mathematics \\ University of Houston \\ Houston, TX 77204-3008 \\USA}
\email{matt\_tobin67@hotmail.com}
 
\address{Department of Mathematics \\ University of Houston \\ Houston, TX 77204-3008 \\USA}
\email{letty@math.uh.edu}
 
\address{Department of Mathematics \\ University of Houston \\ Houston, TX 77204-3008 \\USA}
\email{tomforde@math.uh.edu}

%\thanks{This research was supported by NSF Postdoctoral Fellowship DMS-0201960}
 
\date{\today}
 
 
 
\thanks{Parts of this research were supported by NSF Grant DMS-0807399,
and by DMS-0604429 (REU supplement).}
 
\begin{abstract}
 
We develop the theory of frames and Parseval frames 
for finite-dimen\-sio\-nal
vector spaces over the binary numbers. This includes
characterizations which are similar to frames and Parseval frames 
for real or complex Hilbert spaces, and the discussion of conceptual
differences caused by the lack of a proper inner product on
binary vector spaces. 
We also define switching equivalence for 
binary frames, and list all equivalence classes of binary Parseval frames in lowest 
dimensions, excluding cases of trivial redundancy.
\end{abstract}
 
\maketitle

%%%%%%%%%%%%%%%%%%%%%%%%%%%%%%%%%%%%%%%%%%%%%%%%
\section{Introduction}
%%%%%%%%%%%%%%%%%%%%%%%%%%%%%%%%%%%%%%%%%%%%%%%%

There are many conceptual similarities between frames 
and error-correct\-ing linear codes. Frame theory is concerned with stable linear 
embeddings of Hilbert spaces
obtained from mapping a vector to its frame coefficients \cite{DS,Chr03,HKLW}.
The linear dependencies incorporated in the frame
coefficients of a vector help recover from errors such as noise, quantization and
data loss \cite{GVT98,GKK01,RG03,RG04,PK05}, just as linear codes help recover from 
symbol decoding errors and erasures \cite{MWS77}.
Frame design for specific purposes has been related to optimization
problems of a geometric \cite{CK03,SH03,HP} or even discrete nature \cite{BP05,XZG05,Kal06}, 
including combinatorial considerations that are more commonly
associated with error-correcting codes.
On the other hand, one may ask whether concepts from frame theory yield insights in
the binary setting. This is the motivation of the present paper. 

We translate
many of the essential results on frames for finite-dimensional real or complex Hilbert spaces
to analogous statements for vector spaces over the binary numbers. In the first part,
we show that in the binary case, the spanning property of a family of vectors
is equivalent to having a reconstruction identity with a dual family. This means,
both properties can be used interchangeably as a definition of frames, 
as on finite
dimensional real or complex Hilbert spaces.
On the other hand, we demonstrate that an attempt to define binary frames 
similarly to the real or complex case via norm inequalities
fails in binary vector spaces, because they lack an inner product and a polarization identity. 
In the main part of this paper, we focus on Parseval
frames, which have a particularly simple reconstruction identity. We characterize
binary Parseval frames in terms of their frame operator and develop a notion
of switching equivalence for binary frames, similar to the concept for
real or complex frames \cite{GKK01,HP,BP05}. Moreover, we introduce the notion of trivial redundancy,
caused by repeated vectors or the inclusion of the zero vector in the frame.
Ignoring cases of trivial redundancy and choosing representatives from each switching equivalence class simplifies the enumeration of 
binary Parseval frames. By an exhaustive search, we have found  that if $k \in \{4,5,\dots, 11\}$
then all frames that are not trivially redundant in $\mathbb Z_2^4$ with 
$k$ vectors belong to one switching equivalence class.
Further simplifications for the search of all binary Parseval frames
are obtained from a combinatorial consideration, which might be useful for a future
effort to catalogue binary Parseval frames in larger dimensions. 

The remainder of this paper is organized as follows:
In Section~\ref{sec:prelim}, we define frames for
finite-dimensional binary vector spaces. Section~\ref{sec:Parseval}
specializes the discussion to Parseval frames. Finally,  in Section~\ref{sec:catalog}, we define switching equivalence
for binary frames and give a catalog of representatives from each equivalence class of
Parseval frames in lowest dimensions, excluding the trivially redundant ones.

%%%%%%%%%%%%%%%%%%%%%%%%%%%%%%%%%%%%%%%%%%%%%%%%
\section{Preliminaries}\label{sec:prelim}
%%%%%%%%%%%%%%%%%%%%%%%%%%%%%%%%%%%%%%%%%%%%%%%%

 In this  section we first revisit the essentials of frames over the fields $\R$ or~$\C$, the real 
 or complex numbers.
 We then proceed to develop the concept of frames
 over the field $\mathbb Z_2$, that is, the field with two elements $\{0,1\}$, 
 where $0$ is the neutral element with respect to addition, and $1$ is the neutral
 element with respect to multiplication. The main insight of this section is that 
 while there are equivalent characterizations of certain types of frames when the ground field is 
 $\R$ or $\C$, this is not true over
 $\mathbb Z_2$, because 
 the polarization identity is no longer available due to the lack of an inner product.
 
If $\mathcal H$ is a finite-dimensional Hilbert space over $\R$ or $\C$ with inner product $\langle \cdot, \cdot \rangle$, then a family of vectors $\F := \{ f_1, f_2, \ldots, f_k \}$ in $\mathcal H$ is called a \emph{frame} if there exist real numbers $A$ and $B$ such that $0 < A \leq B < \infty$ and
\begin{equation} \label{frame-condition}
A \| x \|^2 \leq \sum_{j=1}^k | \langle x, f_j\rangle |^2 \leq B \| x \|^2 \qquad \text{ for all $x \in \mathcal H$.}
\end{equation}
The inequalities displayed in \eqref{frame-condition} are known as the \emph{frame condition}, and it can be shown that when $\mathcal H$ is finite dimensional, then  the set $\F$ satisfies the frame condition if and only if $\operatorname{span} \F = \mathcal H$ \cite[Proposition~3.18]{HKLW}.
In this case, there exist vectors $\{g_1, g_2, \ldots, g_k\}$ which provide the
\emph{reconstruction identity}
$$
   x = \sum_{j=1}^k \langle x, f_j \rangle g_j \mbox{ for all } x \in \mathcal H \, .
$$
While the family $\{g_j\}_{j=1}^k$ may not be unique, there is a canonical choice. If
we define the so-called \emph{frame operator} $S$ on $\mathcal H$ by $Sx=\sum_{j=1}^k \langle x, f_j \rangle f_j$,
then setting $g_j = S^{-1} f_j$ for $j \in \{1, 2, \ldots k\}$ yields the reconstruction identity
\cite{Chr03}. The family $\{g_j\}_{j=1}^k$
is also called the \emph{canonical dual frame}.

A frame $\F = \{ f_1, \ldots, f_k \}$ is called a \emph{Parseval frame} (or sometimes a \emph{normalized tight frame}) if we can choose $A = B = 1$ in the frame condition, so that 
\begin{equation} \label{Parseval-Identity}
\sum_{j=1}^k | \langle x, f_j \rangle |^2 = \| x \|^2 \qquad \text{ for all $x \in \mathcal H$.}
\end{equation}
Using the polarization identity, it can be shown (see \cite[Proposition~3.11]{HKLW}) that $\F$ is a Parseval frame if and only if 
\begin{equation} \label{Reconstruction-Formula}
x = \sum_{j=1}^k \langle x, f_j \rangle f_j  \qquad \text{ for all $x \in \mathcal H$.}
\end{equation}
The simple form of the reconstruction formula for Parseval frames
%
%The equation shown in \eqref{Reconstruction-Formula} is called the \emph{reconstruction formula}, and it shows that any vector may be reconstructed as a linear combination of the $f_i$'s using the scalars $\langle x, f_i \rangle$ as coefficients in the linear combination.  This
has many practical uses in engineering and computer science \cite{GVT98,GKK01,K}.
 
 We now turn to frames over the binary numbers.
 
 The first two goals in this paper are to develop the notion of frames and of Parseval frames for finite-dimensional vector spaces over the field $\Z_2$.  Any such vector space has the form $\Z_2^n = \Z_2 \oplus \ldots \oplus \Z_2$ for some $n \in \N$. 
 
 \begin{definition}
  A family of vectors $\mathcal F = \{f_1 , f_2, \ldots f_k\}$ in $\mathbb Z_2^n$ is a \emph{frame}
 if it spans $\mathbb Z_2^n$.
 \end{definition}

 We have chosen this form of the definition because the field $\Z_2$ has no notion of positive elements, so that  it is impossible to find a properly defined inner product, let alone a norm on $\Z_2^n$, which
 would be needed to formulate a direct analogue of the frame condition (\ref{frame-condition}).
 
 Nevertheless, we want to show that an analogue of the reconstruction identity can be deduced
 with the help of a  $\Z_2$-valued ``dot product" in place of an inner product.
 
\begin{definition}
We define a bilinear map $( \cdot, \cdot ) : \Z_2^n \times \Z_2^n \to \Z_2$ called the \emph{dot product} on $\Z_2^n$ by $$\left( \left( \begin{smallmatrix} a_1 \\ \vdots \\ a_n \end{smallmatrix} \right),  \left( \begin{smallmatrix} b_1 \\ \vdots \\ b_n \end{smallmatrix} \right) \right) := \sum_{i=1}^n a_i b_i.$$
\end{definition}
 
We see that the dot product $( \cdot, \cdot )$ is symmetric and $\Z_2$-linear in each component, but it is degenerate: It is possible to have $x \in \Z_2^n$ with $(x,x) = 0$ but $x \neq 0$.  Furthermore, because the dot product is degenerate, it does not provide a norm on $\Z_2^n$.  Nonetheless, we will use the dot product as an analogue of the inner products on $\R^n$ and $\C^n$, and for expressions in $\R^n$ or $\C^n$ involving $\langle x, y \rangle$ or $\| x \|^2$, we shall consider analogous expressions in $\Z_2^n$ involving $(x,y)$ or $(x,x)$, respectively.
 
To establish the equivalence between the spanning property and the reconstruction identity
 for frames, we unfortunately cannot simply use the same strategy as in the real or complex case.
 If we take the dot product instead of an inner product to
 define the frame operator, then the spanning property of the frame does not guarantee
 that the frame operator is invertible. To see this, we note that
 the family $\{1,1\}$ is spanning for $\mathbb Z_2$, but the analogue of the frame operator
 maps every $x\in \mathbb Z_2$ to $x+x=0$. A similar family can be obtained for
 any $\mathbb Z_2^n$, $n\ge 1$, by repeating vectors of an arbitrary spanning set.
 
 To build
 an alternative strategy that relates the spanning property with the existence of a reconstruction identity,
 we first  recall that the dot product mediates a canonical mapping between vectors and linear functionals.

 \begin{lemma}
If $\phi: \Z_2^n\to\Z_2$ is a linear functional then there exists a unique $z\in \Z_2^n $ such 
that $\phi(x) = (x,z)$ for all $x \in\Z_2^n$.
\end{lemma}

\begin{proof}
Let $\phi$ be a linear functional.
Let \{$e_1$,$\dots$,$e_n$\} be the canonical basis for $\Z_2^n$, and
let $z=\phi(e_1)e_1 + \dots+ \phi(e_n)e_n$.
We now observe that if $x\in\Z_2^n$, with $x=\sum_{i=1}^{n}a_i e_i$ for $a_i\in \Z_2$,
then $\phi(x)=\sum_{i=1}^{n}a_i \phi(e_i) = (x, z)\, . $

To verify the uniqueness, assume there is $z'$ such that $\phi(x)=(x,z')=(x,z)$.
Choosing $x$ among the canonical basis vectors gives $\phi(e_i)=(e_i,z')=(e_i,z)$
and thus $z$ and $z'$ are identical.
\end{proof}

 \begin{theorem}
Given a family {$\mathcal F = \{f_j\}_{j=1}^k$} in $\Z_2^n$, then
$\mathcal F$ is a frame if and only if there exist 
vectors $\{g_j\}_{j=1}^k$
%linear functionals {$\gamma_1$, $\gamma_2$, $\dots$, $\gamma_k$} 
such that for all $y$ $\in$  $\Z_2^n$ 
\begin{equation} \label{eq:frame-ident}
  y=\sum_{j=1}^{k} (y,g_j) f_j . 
\end{equation}
 \end{theorem}

\begin{proof}
We note that if (\ref{eq:frame-ident}) is true, then necessarily $\{f_j\}_{j=1}^k$
is spanning.

Conversely, let us assume that $\{f_j\}_{j=1}^k$ is a frame for $\mathbb Z_2^n$.
In a first step, we prove that there are linear functionals $\{\gamma_1$, $\gamma_2$, 
$\dots$, $\gamma_k\}$ such that for all $y \in \Z_2^n$,
$$
   y=\sum_{j=1}^{k} \gamma_j(y) f_j . 
$$
For any family of linear functionals {$\gamma_1$, $\gamma_2$, $\dots$, $\gamma_k$}, we note
that the expression
$\sum_{j=1}^{k} \gamma_j(y)f_j$  
is linear in $y$, so it is enough to show that there exist linear functionals giving
 \begin{center} $w_i$ $=$ $\displaystyle\sum_{j=1}^{k}$ $\gamma_j(w_i) f_j$ 
 for all vectors in some basis  {$w_1$, $\dots$, $w_n$} of $\Z_2^n$   \, .\end{center}
 To establish this,
we choose a subset of {$\{f_1$, $\dots$,$f_k\}$} which is spanning and linearly independent,
that is, a basis. Without loss of generality, assume that this set is
 {$\{f_1$, $\dots$, $f_n\}$}. Choosing the
dual basis {$\{\gamma_1$, $\dots$, $\gamma_n\}$} to {$\{f_1$, $\dots$, $f_n\}$}, 
characterized by 
\begin{center}$\gamma_j(f_i) = \delta_{ij}$ \, , for all $i, j \in \{1, 2, \dots n\}$ \end{center}
we obtain
\begin{center} $\displaystyle\sum_{j=1}^{n}$ $\gamma_j(f_i)f_j$ $=$ $f_i$ \, .\end{center} 
Thus if we enlarge the set $\{\gamma_j\}_{j=1}^n$ by setting $\gamma_j=0$ if $j>n$, 
then
%Thus if we choose $f_j$ $=$ $y_j$ %%need to come back
%then for all $i$ $\epsilon$ {1,$\dots$,$k$ }
$$f_i=\sum_{j=1}^{k}\gamma_j(f_i)f_j $$
 and by linearity
%and because {$f_1$, $\dots$, $f_k$} is a basis for $\Z_2^k$
$$y=\sum_{j=1}^{k}\gamma_j(y)f_j \mbox { for any }y\in \Z_2^n \, .
$$

In the final step of the proof, we apply the preceding lemma which yields
for each $\gamma_j$ a corresponding vector $g_j$ satisfying
$\gamma_j(y) = (y,g_j)$ for all $y \in \mathbb Z_2^n$.
\end{proof}

%%%%%%%%%%%%%%%%%%%%%%%%%%%%%%%%%%%%%%%%%%%%%%%%
\section{Parseval frames for $\Z_2^n$}\label{sec:Parseval}
%%%%%%%%%%%%%%%%%%%%%%%%%%%%%%%%%%%%%%%%%%%%%%%%
 
In this section we present the definition of Parseval frames for $\mathbb Z_2^n$
and illustrate the conceptual differences between such frames in the real or complex case
and in the binary case.
%Because the reconstruction formula (see \eqref{Reconstruction-Formula}) is the most useful formulation of Parseval frames in vector spaces over $\R$ or $\C$, we use this as our definition in $\Z_2^n$.
 
\begin{definition}
A family of vectors $\F = \{ f_1, \ldots, f_k \}$ in $\Z_2^n$ is a \emph{Parseval frame} if
\begin{equation} \label{Reconstruction-Formula-Z2}
x = \sum_{j=1}^k ( x, f_j )f_j  \qquad \text{ for all $x \in \Z_2^n$.}
\end{equation}
\end{definition}
\noindent Observe that a binary Parseval frame necessarily spans $\Z_2^n$, and moreover if $\F$ is a Parseval frame, we must have $k \geq n$.
 
It is natural to ask if, in analogy with the real and complex cases, being a Parseval frame in $\Z_2^n$ is equivalent to having a Parseval identity as in \eqref{Parseval-Identity}.  It turns out that this is not the case.
 
\begin{proposition}
 If $\F= \{ f_1, \ldots, f_k \}$ is a Parseval frame
for $\mathbb Z_2^n$, then 
\begin{equation} \label{Parseval-Identity-Z2}
\sum_{j=1}^k ( x, f_j )^2 = ( x, x) \qquad \text{ for all $x \in \Z_2^n$.}
\end{equation}
However, in general, the converse does not hold.
\end{proposition}
 
\begin{proof}
If $\F$ is a Parseval frame, then using the $\Z_2$-linearity of the first component of the dot product, for any $x \in \Z_2^n$ we have $$(x, x) = \left(  \sum_{j=1}^k (x,f_j)f_j ,x\right) = \sum_{j=1}^k (x,f_j)(f_j,x) = \sum_{j=1}^k (x,f_j)^2.$$  

To see that the converse does not hold in general,  consider $\left( \begin{smallmatrix} 1 \\ 1 \end{smallmatrix} \right) \in \Z_2^2$, then for any $x = \left( \begin{smallmatrix} a_1 \\ a_2 \end{smallmatrix} \right) \in \Z_2^2$ we have $$( x, \left( \begin{smallmatrix} 1 \\ 1 \end{smallmatrix} \right) ) = a_1 + a_2 = a_1^2 + a_2^2 = (x, x). $$  Hence $\F = \left\{ \left( \begin{smallmatrix} 1 \\ 1 \end{smallmatrix} \right) \right\}$ satisfies \eqref{Parseval-Identity-Z2}.  However, $\F$ contains one element, so $\F$ does not span $\Z_2^2$, and $\F$ is not a Parseval frame.
\end{proof}

\begin{remark}
More generally, we can produce counterexamples for any $n\ge 2$, meaning sets
which give the Parseval property without spanning $\mathbb Z_2^n$. First we consider even $n$.
Let $\{f_1, \ldots f_k\}$ be the family of all vectors which contain exactly two
$1$'s. Thus, there are $k=\left({n \atop 2}\right)$ such vectors.
If the fist vector is chosen as $f_1=(1,1,0,\ldots,0)^t$ and
$x=(a_1, a_2, \ldots a_n)^t$, then over $\mathbb Z_2$,
$$
    (x, f_1)^2 = (a_1 + a_2)^2 = a_1^2 + a_2^2 \, .
$$
Evaluating other dot products similarly  gives
$$
   \sum_{j=1}^k (x,f_j)^2 =  \sum_{i=1}^n a_i^2 \, 
$$
because each $a_i^2$ appears in $n-1$ terms in the sum, and $n-1 \mbox{ mod } 2= 1$
by the assumption that $n$ is even.

However, the vectors  $\{f_j\}_{j=1}^k$ are not spanning for $\mathbb Z_2^n$,
because they contain an even number of $1$'s and so
does any linear combination of them. 

If $n$ is odd, then we split $\mathbb Z_2^n = \mathbb Z_2 \oplus \mathbb Z_2^{n-1}$
and construct the above family $\{f_1, f_2 , \ldots f_k\}$ for the second summand. 
Now this family can be enlarged by the first canonical basis vector $e_1$ to 
$\{e_1, f_1, f_2 , \ldots f_k\}$ which has the Parseval property but
is not spanning, because  $\{f_1, f_2 , \ldots f_k\}$ does not span $\mathbb Z_2^{n-1}$.
 \end{remark}

%We can see this by some following counter examples.
 
%\begin{flushleft}
 
%Note that $ "a,b,c,...\in \mathbb{Z}_2$" means $"a,b,c,...\in\{0,1\}$" and the multiple of 2 in $ \mathbb{Z}_2$ is zero.
%\end{flushleft}

%%%%%%%%%%%%%%%%%%%%%%%%%%%%%%%%%%%%%%%%%%%%%%%%
\section{Towards a catalog of binary Parseval frames}\label{sec:catalog}
%%%%%%%%%%%%%%%%%%%%%%%%%%%%%%%%%%%%%%%%%%%%%%%%
 
In principle, all Parseval frames for $\mathbb Z_2^n$ could be catalogued individually, but 
even for relatively small $n$ this is already an extensive list. 
In order to obtain a more efficient
way of enumerating Parseval frames, we use an equivalence relation 
which has been called switching equivalence for real or complex frames \cite{GKK01,HP,BP05}.
It is most easily
formulated in terms of the Grammian of a Parseval frame, as defined below. The catalogue of
frames can then be reduced to representatives of each equivalence class. To prepare the 
definition of the equivalence relation, we discuss certain matrices related to frames.

We write $A \in M_{m,n} (\Z_2)$ when $A$ an $m \times n$ matrix with entries in $\Z_2$. We often view $A$ as a linear map from $\Z_2^m$ to $\Z_2^n$ by left multiplication. In particular, $A\in M_{n}$
denotes an $n\times n$ matrix which is associated with a map from $\mathbb Z_2^n$ to itself.
 We write $A_{i,j}$ for the $(i,j)$\textsuperscript{th} entry of $A$, and we let $A^*$ denote the transpose of $A$; that is, $A^* \in M_{n,m}(\Z_2)$ with $A^*_{i,j} := A_{j,i}$.
By the rules of matrix multiplication, we have $(Ax,y)=(x,A^*y)$ for all $A\in M_n(\Z_2)$.  
 
%\begin{lemma}
%If $A \in M_n(\Z_2)$ and $x,y \in \Z_2^n$, then $(Ax,y) = (x, A^*y)$.
%\end{lemma}
% 
%\begin{proof}
%Write $x= \left( \begin{smallmatrix} a_1 \\ \vdots \\ a_n \end{smallmatrix} \right)$ and $y= \left( \begin{smallmatrix} b_1 \\ \vdots \\ b_n \end{smallmatrix} \right)$.  Then 
%\begin{align*}
%(Ax,&y) = \left(  \left( \begin{smallmatrix} \sum_{j=1}^n A_{1,j} a_j \\ \vdots \\ \sum_{j=1}^n A_{n,j} a_j \end{smallmatrix} \right),  \left( \begin{smallmatrix} b_1 \\ \vdots \\ b_n \end{smallmatrix} \right) \right) = \sum_{i=1}^n \sum_{j=1}^n A_{i,j} a_j b_i =  \sum_{j=1}^n \sum_{i=1}^n a_j A_{i,j} b_i \\
%& =  \left(  \left( \begin{smallmatrix} a_1 \\ \vdots \\ a_n \end{smallmatrix} \right), \left( \begin{smallmatrix} \sum_{i=1}^n A_{i,1}^* b_i \\ \vdots \\ \sum_{i=1}^n A_{j,n}^* b_i \end{smallmatrix} \right)
%\right) =  \left(  \left( \begin{smallmatrix} a_1 \\ \vdots \\ a_n \end{smallmatrix} \right), \left( \begin{smallmatrix} \sum_{i=1}^n A^*_{1,i} b_i \\ \vdots \\ \sum_{i=1}^n A^*_{n,j} b_i \end{smallmatrix} \right)
%\right) = (x, A^* y).
%\end{align*}
%\end{proof}
 
\begin{definition}
If $U \in M_n(\Z_2)$, then we say $U$ is a \emph{unitary} if $U$ is invertible and $U^{-1} = U^*$.
\end{definition}

\begin{lemma} \label{dot-prod-lem}
If $x \in \Z_2^n$ and $(x,y) = 0$ for all $y \in \Z_2^n$, then $y=0$.
\end{lemma}
 
\begin{proof}
Write $x= \left( \begin{smallmatrix} a_1 \\ \vdots \\ a_n \end{smallmatrix} \right)$.  If $\{e_1, \ldots, e_n \}$ is the standard basis for $\Z_2^n$, then for all $1 \leq i \leq n$ we have $a_i = ( x, e_i ) = 0$.  Thus $x = 0$.
\end{proof}

\begin{proposition} \label{unitary-preserve-dot-prod-prop}
Let $U \in M_n(\Z_2)$, then $U$ is a unitary if and only if for all $x,y \in \Z_2^n$ we have $(Ux, Uy) = (x,y)$.
\end{proposition}
 
\begin{proof}
If $U$ is a unitary, then $U^* = U^{-1}$ and for all $x,y \in \Z_2^n$ we have $(Ux,Uy) = (x, U^*Uy) = (x,Iy) = (x,y)$.  Conversely, if for all $x, y \in \Z_2^n$ we have $(Ux, Uy) = (x,y)$, then for a given $x \in \Z_2^n$ we see that $(U^*Ux,y) = (Ux,Uy) = (x,y)$ for all $y \in \Z_2^n$, and Lemma~\ref{dot-prod-lem} implies that $U^*Ux = x$.  Since $x$ was arbitrary, this shows that $U^*U = I$, and because $U$ is square, we have that $U$ is invertible and $U^{-1} = U^*$.
\end{proof}

%\begin{remark}
In contrast to unitaries for Hilbert spaces over
$\mathbb F = \R \text{ or } \C$, the condition $\langle Ux,Ux \rangle = \langle x, x \rangle$ for all $x \in \mathbb F^n$ is not equivalent to unitarity when the field $\mathbb F$ is $\mathbb Z_2$.
%let $U \in M_n(F)$, and let $U^*$ denotes the conjugate transpose of $U$.  Then $U$ is invertible with $U^{-1} = U^*$ if and only if .  (This is a consequence of the polarization identity.) However, an analogous statement does not hold for $\Z_2^n$. 
%\end{remark}

We have the following counterexamples for $n\ge 2$

\begin{proposition}
For any $n \ge 2$, there exist $A \in M_n(\mathbb Z_2)$ such that 
$ (Ax,Ax)=(x,x)$ for all $x \in \mathbb Z_2^n$
but $A$ is not invertible, and thus not unitary.
\end{proposition}
\begin{proof}
We define the matrix $A$ by 
$$
    A_{i,j} = \begin{cases} 1, & \mbox{ if } i=j=1 \mbox{ or } j-i=1\, ,\\
                                             0, & \mbox{ else. }
                                             \end{cases}
$$
This means, the last row of $A$ contains only zeros and thus
$A$ does not have full rank and is not invertible.

However, given $x=(a_1, a_2, \ldots a_n)^t$
we have 
$$ 
 A x =  \left(\begin{smallmatrix}a_1 + a_2\\ a_3\\ a_4\\ \dots \\ a_n\\ 0\end{smallmatrix}\right)   
$$ 
and thus
$$
  (Ax,Ax) = (a_1+a_2)^2 + a_3^2 + \dots + a_n^2
    = \sum_{i=1}^n a_i^2 =(x,x) \, .
$$
\end{proof}
 
\begin{definition}
Let $\F = \{ f_1, \ldots, f_k \} \subseteq \Z_2^n$.  The \emph{analysis operator} for $\F$ is the $k \times n$ matrix containing the frame vectors as rows,
$$\Theta_\F = \begin{pmatrix} \leftarrow & f_1 & \rightarrow \\ & \vdots & \\  \leftarrow & f_k & \rightarrow \end{pmatrix}.$$ The \emph{synthesis operator} for $\F$ is the $n \times k$ matrix $$\Theta^*_\F = \begin{pmatrix} \uparrow &  & \uparrow \\ f_1 & \cdots & f_k \\  \downarrow & & \downarrow \end{pmatrix},$$ with the elements of $\F$ as columns.  The \emph{frame operator} for $\F$ is the $n \times n$ matrix $$S_\F := \Theta_\F^* \Theta_\F,$$ and the \emph{Grammian operator} for $\F$ is the $k \times k$ matrix $$G_\F :=  \Theta_\F \Theta_\F^*.$$  Note that $(G_\F)_{i,j} = (f_j, f_i)$ for all $1 \leq i,j \leq k$.  When there is no ambiguity in the choice of $\mathcal F$, we shall omit the $\F$ subscript on these matrices and simply write $\Theta$, $\Theta^*$, $S$, and $G$.
\end{definition}
 
\begin{theorem} \label{Parseval-analysis-op-thm}
Let $\F = \{ f_1, \ldots, f_k \} \subseteq \Z_2^n$, then $\F$ is a Parseval frame if and only if $S_\F$ is equal to the identity matrix.
\end{theorem}
 
\begin{proof}
Let $\{e_1, \ldots, e_n \}$ be the standard basis for $\Z_2^k$.  
Observe that for any $x \in \Z_2^n$ we have $\Theta_\F x = \sum_{i=1}^k (x, f_i) e_i$.  Also, for any $1 \leq i \leq n$ we have $\Theta_\F^* e_i = f_i$.  Thus we have $$S_\F x = \Theta_\F^* \Theta_\F x = \Theta_\F^*\left( \sum_{i=1}^k (x, f_i) e_i \right) = \sum_{i=1}^k (x, f_i) \Theta_\F^* e_i =  \sum_{i=1}^k (x, f_i) f_i.$$  It follows that $\sum_{i=1}^k (x, f_i) f_i = x$ for all $x \in \Z_2^n$ if and only if $S_\F x = x$ for all $x \in \Z_2^n$.  Thus $\F$ is a Parseval frame if and only if $S_\F$ is the identity matrix.
\end{proof}

\begin{definition}
Given two families $\F = \{f_1, \ldots, f_k \}$ and $\G = \{ g_1, \ldots, g_n \}$ in $ \Z_2^n$, then we say \emph{$\F$ is unitarily equivalent to $\G$} if there exists a unitary $U \in M_n(\Z_2)$ such that $Uf_i = g_i$ for all $1 \leq i \leq k$.
\end{definition}
 
It is easy to show that unitary equivalence is an equivalence relation.
 
\begin{proposition}
Let $\F = \{f_1, \ldots, f_k \} \subseteq \Z_2^n$ and $\G = \{ g_1, \ldots, g_n \} \subseteq \Z_2^n$ be Parseval frames, then $\F$ is unitarily equivalent to $\G$ if and only if $G_\F = G_\G$.
\end{proposition}
 
\begin{proof}
Since $\F$ and $\G$ are Parseval frames, it follows from Theorem~\ref{Parseval-analysis-op-thm} that $S_\F$ and $S_\G$ are the identity matrices.  Suppose that $G_\F = G_\G$.  Define $U$ to be the $n \times n$ matrix $U := \Theta_\G^* \Theta_\F$, then 
\begin{align*}
U^*U &= (\Theta_\G^* \Theta_\F)^* \Theta_\G^* \Theta_\F = \Theta_\F^* \Theta_\G \Theta_\G^* \Theta_\F = \Theta_\F^* G_\G \Theta_\F \\ 
&= \Theta_\F^* G_\F \Theta_\F =  \Theta_\F^* \Theta_\F \Theta_\F^* \Theta_\F = S_\F S_\F = I.
\end{align*}
Since $U$ is square, it follows that $U$ is invertible and $U^{-1} = U^*$, so that $U$ is a unitary.  Furthermore,
$$U \Theta_\F^* =  \Theta_\G^* \Theta_\F \Theta_\F^* = \Theta_\G^* G_\F = \Theta_\G G_\G = \Theta_\G^* \Theta_\G \Theta_\G^* = S_\G \Theta_\G^* = \Theta_\G^*.$$  Thus $U$ times the $i$\textsuperscript{th} column of $\Theta_\F^*$ is equal to the $i$\textsuperscript{th} column of $\Theta_\G^*$.  Thus for all $1 \leq i \leq k$ we have $U f_i = g_i$, so that $\F$ and $\G$ are unitarily equivalent.
 
Conversely, if $\F$ and $\G$ are unitarily equivalent, then there exists a unitary $U \in M_n(\Z_2)$ such that $U f_i = g_i$ for all $1 \leq i \leq k$.  Thus Proposition~\ref{unitary-preserve-dot-prod-prop} implies that $$(G_\F)_{i,j} = (f_j, f_i) = (Uf_j, Uf_i) = (g_j, g_i) = (G_\G)_{i,j}.$$  Hence $G_\F = G_\G$.
\end{proof}

\begin{example}
We present two examples of unitary equivalence:
\begin{align*}
\F= \left\{ \left( \begin{smallmatrix} 0\\1\\0\\0\\1 \end{smallmatrix} \right),
\left( \begin{smallmatrix} 0\\1\\0\\1\\1 \end{smallmatrix} \right),
\left( \begin{smallmatrix} 0\\1\\1\\0\\1 \end{smallmatrix} \right),
\left( \begin{smallmatrix} 1\\0\\1\\1\\1 \end{smallmatrix} \right),
\left( \begin{smallmatrix} 1\\1\\0\\0\\1 \end{smallmatrix} \right),
\left( \begin{smallmatrix} 1\\1\\1\\1\\0 \end{smallmatrix} \right) \right\}
\end{align*}
\begin{align*}
\mathcal H= \left\{ \left( \begin{smallmatrix} 0\\1\\0\\1\\0 \end{smallmatrix} \right),
\left( \begin{smallmatrix} 0\\1\\0\\1\\1 \end{smallmatrix} \right),
\left( \begin{smallmatrix} 0\\1\\1\\1\\0 \end{smallmatrix} \right),
\left( \begin{smallmatrix} 1\\0\\1\\1\\1 \end{smallmatrix} \right),
\left( \begin{smallmatrix} 1\\1\\0\\1\\0 \end{smallmatrix} \right),
\left( \begin{smallmatrix} 1\\1\\1\\0\\1 \end{smallmatrix} \right) \right\}
\end{align*}
Computing the Grammian operator $G$ for both $\F$ and $\mathcal H$ we find:
\begin{align*}
G_\F = G_{\mathcal H} = \begin{bmatrix} 0&0&0&1&0&1\\0&1&0&0&0&0\\0&0&1&0&0&0\\1&0&0&0&0&1\\0&0&0&0&1&0\\1&0&0&1&0&0\end{bmatrix}.
\end{align*}
First notice the structure of $\Theta^*$ created by both $\F$ and $\mathcal H$.
$$\Theta_\F^*=\begin{bmatrix} 0&0&0&1&1&1\\1&1&1&0&1&1\\0&0&1&1&0&1\\0&1&0&1&0&1\\1&1&1&1&1&0\end{bmatrix}
\Theta_{\mathcal H}^*=\begin{bmatrix} 0&0&0&1&1&1\\1&1&1&0&1&1\\0&0&1&1&0&1\\1&1&1&1&1&0\\0&1&0&1&0&1\end{bmatrix}$$
The fourth and fifth rows have swapped places, so naturally one would expect the unitary operator to reflect that.  In fact, the proof gives a direct way to compute $U$, namely if $f_i= U h_i$ then $U=\Theta_\F^* \Theta_{\mathcal H}$. Computing $U$ confirms this.
$$ \Theta_\mathcal{F}^* \Theta_\mathcal{H}=
\begin{bmatrix}1&0&0&0&0\\0&1&0&0&0\\0&0&1&0&0\\0&0&0&0&1\\0&0&0&1&0 \end{bmatrix}$$

Looking at one last example, take $\F$ and $\mathcal H$, two Parseval Frames found in $\Z_2^5$ with six elements:
\begin{align*}
\F= \left\{
\left( \begin{smallmatrix} 0\\0\\0\\1\\1 \end{smallmatrix} \right),
\left( \begin{smallmatrix} 0\\0\\1\\0\\1 \end{smallmatrix} \right),
\left( \begin{smallmatrix} 0\\1\\0\\0\\1 \end{smallmatrix} \right),
\left( \begin{smallmatrix} 1\\0\\0\\0\\1 \end{smallmatrix} \right),
\left( \begin{smallmatrix} 1\\1\\1\\1\\0 \end{smallmatrix} \right),
\left( \begin{smallmatrix} 1\\1\\1\\1\\1 \end{smallmatrix} \right)\right\}
\end{align*}
\begin{align*}
\mathcal H=\left\{
\left( \begin{smallmatrix} 0\\1\\1\\1\\1 \end{smallmatrix} \right),
\left( \begin{smallmatrix} 1\\0\\0\\0\\1 \end{smallmatrix} \right),
\left( \begin{smallmatrix} 1\\0\\0\\1\\0 \end{smallmatrix} \right),
\left( \begin{smallmatrix} 1\\0\\1\\0\\0 \end{smallmatrix} \right),
\left( \begin{smallmatrix} 1\\1\\0\\0\\0 \end{smallmatrix} \right),
\left( \begin{smallmatrix} 1\\1\\1\\1\\1 \end{smallmatrix} \right)\right\}
\end{align*}
Here, while not quite as obvious, differences in the two Parseval frames can be expressed in terms of row manipulations of the synthesis operator which amount to left multiplication  with 
a unitary $U$, $\Theta_\F^*=U\Theta_{\mathcal H}^*$.
\end{example} 

We introduce an additional way to identify frames which coarsens the equivalence relation.

\begin{definition}
Two families $\F=\{f_1, f_2, \dots, f_k\}$ and $\mathcal G = \{g_1, g_2, \dots g_k\}$ in $\mathbb Z^n$
are called \emph{switching equivalent} if there is a unitary $U$ and a permutation $\pi$ of the set
$\{1, 2, \dots k\}$ such that
$$
   f_j = U g_{\pi(j)} \mbox{ for all } j \in \{1, 2, \dots \} \, .
$$ 
\end{definition}

\begin{theorem}
Two Parseval frames $\mathcal F$ and $\mathcal H$ are switching equivalent
if and only if there exists a permutation $\pi$ of the index set such that
$(G_\F)_{i,j} = (G_{\mathcal H})_{\pi(i),\pi(j)}$.
\end{theorem}
\begin{proof}
We note that the condition on the Grammians amounts to
the identity
$$
   G_\F = M G_{\mathcal H} M^*
$$
for a permutation matrix with entries
$$
   M_{i,j} = \begin{cases} 1, & \mbox{ if } \pi(i)=j \\ 0, & \mbox{ else }.\end{cases} 
$$
Being identical up to conjugation by permutation matrices defines an
equivalence relation for Grammians, and thus for frames, which is coarser than unitary equivalence.

Moreover, with a similar proof as in the preceding proposition, we see that
the two Grammians are related by conjugation with a permutation matrix $M$ if and only if
there exists a unitary $U$ such that
$$
   \Theta_\F^* = U \Theta_{\mathcal H}^* M^* \, .
$$
\end{proof}
 
 Apart from switching equivalence, there are other simple ways in which two Parseval frames
 can be related to each other. For example, adding zero vectors to a Parseval frame
 gives another Parseval frame. Moreover, adding pairs of arbitrary vectors to a Parseval
 frame preserves the Parseval property. In both cases, we have artificially increased the 
 redundancy by enlarging the frame.
 In our catalog, we discard Parseval frames with such a trivial source of redundancy.

\begin{definition}
A Parseval frame $\mathcal F=\{f_1, f_2, \dots, f_k\}$ for $\mathbb Z_2^n$ is called \emph{trivially redundant}
if there is $j \in \{1, 2, \dots k\}$ with $f_j=0$, or if there are two indices $i\neq j$
with $f_i=f_j$.
\end{definition}
%taking the complement of a frame in the set of all vectors in $\mathbb Z_2^n$.

After repeated vectors are removed,
Parseval frames can be interpreted as sets of vectors. We consider the set-theoretic complement
of such a Parseval frame.

\begin{theorem}
Let $n\ge 3$.
Let $\mathcal{F}=\left\{f_i\right\}_{i=1}^k$ be a family without repeated vectors in $\mathbb Z_2^n$ 
 and $\mathcal{G}=\mathbb Z_2^n\setminus \mathcal F$.  If $\mathcal{F}$ is a Parseval Frame for $\Z_2^n$, then $\mathcal G$ is also a Parseval frame.
\end{theorem}
\begin{proof}
Let $\mathcal X=\left\{x\in\Z_2^n\right\}$, then we count $2^n-1$ nonzero elements in 
$\mathcal X$.  Thinking of $\mathcal X$ as a sequence $\left\{f_i\right\}_{i=0}^m$ where $m=2^n-1$
and the entries of $f_i$ are given by the binary expansion of $i$, let 
$\Theta_{\mathcal X}^*$ be the matrix with $f_i$ as the $i^{th}$ column.

By simple counting, there are $2^{n-1}$ elements with $1$ in the $i^{th}$ position. 
This means, in each row of $\Theta_{\mathcal X}^*$ the number $1$ appears exactly $2^{n-1}$ times. 
Furthermore there are $2^{n-2}$ elements with $1$ in the $i^{th}$ and $j^{th}$ position, similarly
for any fixed choice of $1$ or $0$ in the  $i^{th}$ and $j^{th}$ position. If $n\ge 3$, then
$2^{n-2}$ is even and consequently,
the dot product of any row of $\Theta_{\mathcal X}^*$ with itself or any other row is equal to $0$, i.e. 
$\Theta_{\mathcal X}^* \Theta_{\mathcal X} = 0$.  

If $\mathcal{F}$ is a Parseval frame, then $\Theta_\F^*\Theta_\F = I$ which implies via the matrix product 
that there is an odd number of elements in $\mathcal{F}$ with $1$ in the $i^{th}$ position, and that among the elements with a $1$ in the $i^{th}$ position there is an even number of elements with a $1$ in the $j^{th}$ position.  

As remarked above, in the entire space there is an even number of elements with $1$ in the $i^{th}$ postion and an even number of elements with $1$ in the $i^{th}$ and $j^{th}$ position.  Thus there is an odd number of elements in the complement $\mathcal G=\mathcal X \setminus\mathcal F$ with $1$ in the $i^{th}$ position and an even number of such elements with $1$ in the $j^{th}$ position, that is $\Theta_\G^*\Theta_\G = I$.  Hence $\G$ is a Parseval Frame.
\end{proof}

\begin{corollary} \label{cor:comp}
If $\mathcal F$ is a Parseval frame for $\mathbb Z_2^n$ 
which is not trivially redundant,
and $\mathcal G$ 
is its set-theoretic complement, then both
$\mathcal F$ and $\mathcal G\setminus \{0\}$ are Parseval frames and
one of them contains at most $2^{n-1}-1$ vectors.
\end{corollary} 
\begin{proof}
After removing the zero vector from $\mathcal G$, the union of both Parseval frames 
$\mathcal F$ and $\mathcal G\setminus\{0\}$ contains
$2^{n}-1$ vectors. This implies that
one of the two frames contains at most half this number, meaning at most $2^{n-1}-1$ vectors.
\end{proof}

To complete the catalog of binary Parseval frames for $\mathbb Z_2^n$, it is only necessary to 
find one representative from
each switching equivalence class of Parseval frames with at most $2^{n-1}-1$ vectors.
Once these Parseval frames have been found, their complements complete the list, because
the switching equivalence of a pair frames 
is equivalent to that of their complements.

\begin{proposition}\label{prop:comp-equiv}
Given
two frames that are not trivially redundant, then they are switching equivalent if and only if 
their set-theoretic complements are.
\end{proposition}
\begin{proof}
This is a consequence of the fact that unitaries are one-to-one maps
on the set $\mathbb Z_2^n$. Thus, if a unitary $U$ maps a frame $\mathcal F$ to a frame $\mathcal G$,
then it also maps the complement of $\mathcal F$ to the complement of $\mathcal G$,
and vice versa.
\end{proof}

We conclude with Table~\ref{tab:one}, a complete list of representatives of switching-equi\-va\-len\-ce classes
of Parseval frames for $n=3$ and $n=4$, excluding ones that are trivially redundant.
Each frame vector in our list is recorded  by the integer obtained from the binary expansion
with the entries of the vector. For example, if a frame vector in $\mathbb Z_2^4$ is $f_1=(1,0,1,1)$,
then it is represented by the integer $2^0+2^2+2^3=13$. Accordingly,
in $\mathbb Z_2^4$, the standard basis is recorded as the sequence of numbers
$1,2,4,8$. 

As explained above, the part of the table with $k>2^{n-1}-1$ vectors,
has been obtained by taking complements of Parseval frames with $k\le 2^{n-1}-1$
vectors, according to Corollary~\ref{cor:comp} and Proposition~\ref{prop:comp-equiv}.
An exhaustive search shows that there is only one switching equivalence class for 
$n=3$ and $k\in \{3,4\}$ and for $n=4$ and
each $k\in\{4,5,6,7\}$, consequently also for
$k \in \{8,9,10,11\}$.

\begin{table}[h]
\begin{tabular}{|r|r|cccc|}
\hline $n$ & $k$ &\multicolumn{4}{c|}{vectors}\\ \hline
3 &  3 & 1 &2 & 4 &\\ \hline
 & 4 & 3 & 5 & 6 & 7 \\ \hline
\end{tabular}
%\end{table}

\vskip1em
%\begin{table}[h]
\begin{tabular}{|r|r|ccccccccccc|}
\hline $n$ & $k$ & \multicolumn{11}{c|}{vectors}\\ \hline
4 & 4 & 1 & 2 & 4 & 8 & & & & & & &\\ \hline
  &  5 & 1 &6 & 10 & 12 & 14& & & & & &\\ \hline
& 6 & 1 & 3 & 5 & 9 & 14 & 15 & & & & & \\ \hline
& 7 & 1 & 2 & 3 & 7 & 11 & 12 & 15 & & & & \\ \hline
 & 8 & 4& 5& 6& 8& 9& 10& 13& 14& & & \\ \hline
 & 9 & 2& 4& 6& 7& 8& 10& 11& 12& 13& & \\ \hline
 & 10 & 2& 3& 4& 5& 7& 8& 9& 11& 13& 15& \\ \hline
 & 11 & 3& 5& 6& 7&9 &10 &11 &12 &13 &14 &15 \\ \hline
\end{tabular}
\vskip1em
\begin{caption}{
Representatives of all switching-equivalence classes, excluding trivially redundant
Parseval frames, for $\mathbb Z_2^3$ and $\mathbb Z_2^4$. 
%Each 
%vector is represented by the decimal which has the entries as digits in its binary expansion.
\label{tab:one}}
\end{caption}

\end{table}

\end{document}